\documentclass[11pt]{article}
\usepackage{amstext}
\usepackage{amssymb}
\usepackage{latexsym}
\newcommand\RR{{\mathbb R}}
\newcommand\CC{{\mathbb C}}
\newcommand\QQ{{\mathbb Q}}

\newcommand\ZZ{{\mathbb Z}}
\newcommand\TT{{\mathbb T}}

\def\d={\,:=\,}

\newcommand{\semdir}
{\rtimes}

\font\frakten=eufm10
\newfam\frakfam
\textfont\frakfam=\frakten

\newtheorem{thm}{Theorem}
\newtheorem{lemma}[thm]{Lemma}
\newtheorem{cor}[thm]{Corollary}
\newtheorem{prop}[thm]{Proposition}
\newtheorem{Defn}[thm]{Definition}
\newtheorem{Ex}[thm]{Example}
\newtheorem{Rem}[thm]{Remark}
\newtheorem{Exs}[thm]{Examples}
\newtheorem{Rems}[thm]{Remarks}
\newtheorem{Defrem}[thm]{Definition and Remark}
\newtheorem{Remnt}[thm]{}
\newenvironment{defn}
 {\begin{Defn} \begin{rm}} {\end{rm} \hfill $\Box$ \end{Defn}}

\newenvironment{ex}
 {\begin{Ex} \begin{rm}} {\end{rm} \hfill $\Box$ \end{Ex}}

\newenvironment{rem}
 {\begin{Rem} \begin{rm}} {\end{rm} \hfill $\Box$ \end{Rem}}

\newenvironment{prf} {{\bf Proof.}}{\hfill $\Box$}

\begin{document}

\title{The weak Paley-Wiener property for group extensions}
\author{Hartmut F\"uhr \\ Institute of Biomathematics and Biometry \\
 GSF Research Center for Environment and Health \\D--85764 Neuherberg \\
 email : fuehr@gsf.de
 }
\date{\today}
\maketitle

\begin{abstract} The paper studies weak Paley-Wiener properties for 
group extensions by use of Mackey's theory. The main theorem establishes 
sufficient conditions on the dual action 
to ensure that the group has the weak Paley-Wiener property.
The theorem applies to yield the weak Paley-Wiener property for large 
classes of simply connected, connected solvable Lie groups (including 
exponential Lie groups), but also criteria for non-unimodular 
groups or motion groups.
\end{abstract}

\setcounter{section}{-1}
\section{Introduction}

The {\em weak Paley-Wiener property (wPW)} can be formulated
as follows: Let $G$ be a second countable, type I locally compact group, and define
\[ {\rm L}^\infty_c(G) = \{ f: G \to \CC : 
f \mbox{ measurable, bounded and compactly supported} \} ~~.\]
Then $G$ has wPW if 
\begin{equation} \label{eqn:wPW_reals}
 \forall \varphi \in {\rm L}_c^\infty(G) ~ \forall 
 \Gamma \subset \widehat{G} :\left( \widehat{f}|_{\Gamma} = 0 \mbox{ and }
 \nu_G(\Gamma) > 0 \Rightarrow f=0 \right)
 \end{equation}
where $\widehat{f}$ denotes the operator-valued Fourier transform
on ${\rm L}^1(G)$, and the measure $\nu_G$ on $\widehat{G}$ is the 
Plancherel measure of $G$.

The statement originated from real Fourier analysis. The fact that
$\RR$ has wPW follows easily from the observation that the Fourier 
transform of a compactly supported function is analytic. An even simpler
argument works for $\ZZ$; Fourier transforms of elements of 
${\rm L}_c^\infty(\ZZ)$ are trigonometric polynomials. On the other hand,
the circle group $\TT$ 
does not have wPW, for obvious reasons.
wPW has been proved for various (unimodular, type I) 
groups, in particular connected, simply connected nilpotent 
Lie groups \cite{Mo,Pa,Gar1,LipRos,ArnLud}; with generalisations to 
completey solvable groups \cite{Gar2}.

The property can be viewed as an uncertainty principle: If
$f$ is compactly supported, $\widehat{f}$ is spread over all of
$\widehat{G}$. It is interesting to compare wPW to
the so-called {\bf qualitative uncertainty
property (Q.U.P.)} studied in \cite{Hog,EchKanKum,
ArnLud}, stating for all $f \in {\rm L}^1(G)$ that
\begin{equation} \label{eqn:QUP}
 \mu_G({\rm supp}(f))  + \nu_G({\rm supp}(\widehat{f})) < \infty
 \Rightarrow f = 0~~,
\end{equation}
where $\mu_G$ is left Haar measure, $\nu_G$ is Plancherel measure
and ${\rm supp}$ denotes the measure-theoretic support 
(unique up to a null set). By contrast to wPW,
this property is obviously not invariant under choice of equivalent 
Plancherel measures, but rather rests on a canonical choice (which
is available for unimodular groups). 

Throughout this paper $G$ denotes a type I locally compact group and $N$ 
a closed normal subgroup of type I, which is in addition regularly embedded. 
All groups are assumed to be second countable.
The aim is to give sufficient criteria
that $G$ has wPW, in terms of analogous properties of $N$
and the little fixed groups.

The paper is structured as follows: We first give a review
of the Mackey machine and its uses for the computation of Plancherel
measure via techniques due to Kleppner and Lipsman. We then present
explicit formulas for the induced representations arising in the 
construction of $\widehat{G}$, and for the associated representations 
of ${\rm L}^1(G)$, which act via certain operator-valued
integral kernels. 
These formulas will allow to prove the main result of this paper,
Theorem \ref{thm:main}, essentially by repeated application of Fubini's 
theorem. 
In the final section we apply Theorem \ref{thm:main} to prove
wPW for a large class of simply connected, connected solvable Lie groups
(Theorem \ref{thm:wPW_solv}),
thereby considerably generalising the previously published
results for nilpotent and completely solvable groups given in
\cite{Mo,Pa,Gar1,LipRos,ArnLud,Gar2}.
Further consequences are criteria for
nonunimodular groups (Corollary \ref{cor:nonunimod}), and a characterisation 
of wPW for motion groups (Theorem \ref{thm:motion}).

\section{Plancherel measure of group extensions}

We follow the exposition in \cite{Li}. For further details and notation
not explained here, 
the reader is referred to this monograph. Throughout the
paper we assume that $G$ is a type I group,
and that $N \lhd G$ is a regularly embedded, type I normal subgroup.
Left Haar measure on a locally compact group is denoted by $| \cdot |$, and
integration against left Haar measure by $\int_G \cdot ~dx$.  
Given a multiplier $\omega$ on $G$, $\widehat{G}^\omega$ denotes the
(equivalence classes of) irreducible unitary $\omega$-representations;
$\omega$ is omitted when it is trivial. $\overline{\omega}$ is the
multiplier obtained by complex conjugation. A multiplier $\omega$ is called 
type I if all $\omega$-representations generate type I von Neumann algebras.

Since we are also dealing with nontrivial Mackey-obstructions, we 
need to recall a multiplier version of the Plancherel theorem.
Given a multiplier $\omega$ on $G$, we denote by $\lambda_{G,\omega}$
the $\omega$-representation of $G$, acting on ${\rm L}^2(G)$ via
\[
 \lambda_{G,\omega} (x) f(y) = \omega(y^{-1},x) f(x^{-1}y) ~~.
\]
If $\omega$ is type I, there exists a 
{\bf Plancherel measure} (unique up to equivalence) $\nu_{G,\omega}$ on
$\widehat{G}^\omega$ decomposing $\lambda_{G,\omega}$,
\[
 \lambda_{G,\omega} \simeq \int_{\widehat{G}^\omega}^\oplus 
 {\rm dim}(\pi) \cdot \pi ~d\nu_{G,\omega} (\pi) 
\]
The associated Fourier transform is given by
\[ {\mathcal F}_\omega : {\rm L}^1(G) \ni f \mapsto (\rho(f))_{\rho
 \in \widehat{G}^\omega} ~~.\]
$\omega$ is omitted whenever it is trivial.

For completeness, we mention the construction of the Plancherel transform
associated to the Plancherel measure. In the unimodular case, the measure
$\nu_G$ can be chosen so that for $f \in {\rm L}^1(G) \cap {\rm L}^2(G)$
the operator field ${\mathcal F}_\omega(f)$ is in fact a field of Hilbert-Schmidt
operators satisfying
\[
 \int_{\widehat{G}} \| {\mathcal F}_\omega(f)(\sigma) \|_{HS}^2 d\nu_G(\sigma) = \|
 f \|_2^2~~,
\]
and the Fourier transform extends by density to a unitary equivalence between
${\rm L}^2(G)$ and the direct integral of Hilbert-Schmidt spaces. The
nonunimodular setting requires right multiplication of
${\mathcal F}_\omega(f)$ with a field of unbounded, densely defined 
selfadjoint operators with densely defined inverse. In particular, 
it does not matter
whether we formulate wPW with reference to operator-valued Fourier- or 
Plancherel transform, though the first one is obviously simpler.

Mackey's theory rests on the dual action of $G$ on $\widehat{N}$.
We assume that $N$ is regularly embedded, i.e., the orbit space 
$\widehat{N}/G$ is countably separated. For $\gamma \in \widehat{N}$ let 
$G_\gamma$ denote the fixed group under the dual action. Now
Mackey's theory provides $\widehat{G}$ as the disjoint union
\[
 \bigcup_{G.\gamma \in \widehat{N}/G} \{ {\rm Ind}_{G_\gamma}^G
 \gamma' \otimes \rho : \rho \in \widehat{G_\gamma/N}^{\overline{\omega_\gamma}} \}~~.
\]
Here $\omega_\gamma $ denotes the multiplier associated to $\gamma$, and
$\gamma'$ denotes a fixed choice of an $\omega_\gamma$-representation of 
$G_\gamma$ on ${\mathcal H}_\gamma$ extending $\gamma$.
Finally, $\rho \in \widehat{G_\gamma/N}^{\overline{\omega_\gamma}}$ is
identified with its 
``lift'' to $G_\gamma$. In the following, we use the notation
\[
 \pi_{G. \gamma, \rho} =  {\rm Ind}_{G_\gamma}^G
 \left( \gamma' \otimes \rho \right) 
\] 
Note that under our assumptions, it follows from a theorem due to
Mackey that all 
$(G_\gamma/N,\overline{\omega_\gamma})$ are type I (see e.g.
\cite[Chapter III, Theorem 4]{Li}).
Thus we have description of $\widehat{G}$ 
as a ``fibred set'', with base space $\widehat{N}/G$ and fibre
\[
 \{ {\rm Ind}_{G_\gamma}^G
 \gamma' \otimes \rho : \rho \in \widehat{G_\gamma/N}^{\overline{\omega_\gamma}} \}
 \simeq \widehat{G_\gamma/N}^{\overline{\omega_\gamma}}
\]
associated to $G.\gamma$. Moreover, there exist Plancherel measures
on $\widehat{N}$ as well as on 
$\widehat{G_\gamma/N}^{\overline{\omega_\gamma}}$. 
Now $\nu_G$ is obtained by taking the projective Plancherel measures 
in the fibres,
and ``glueing'' them together using a pseudo-image $\overline{\nu}$
of $\nu_N$ on $\widehat{N}$, i.e. the quotient measure of a finite
measure equivalent to $\nu_N$. In formulas, up to equivalence $\nu_G$
is given according to \cite[I, 10.2]{KlLi} by 
\begin{equation} \label{eqn:Pl_meas}
 d\nu_G(\pi_{G.\gamma,\rho}) = d\nu_{\widehat{{G_\gamma}/N}^{\overline{\omega_\gamma}}} (\rho)
 d\overline{\nu}(G. \gamma)~~.
\end{equation} 
But $\overline{\nu}$ also enters in a measure decomposition
of $\nu_N$: There exists quasi-invariant measures
$\mu_{G.\gamma}$ on the orbits $G.\gamma \in \widehat{N}/G$ such that 
\begin{equation} \label{eqn:meas_decomp_Pl}
 d\nu_N (\pi) = d\mu_{G.\gamma} (\pi) d\overline{\nu}(G.\gamma) ~~.
\end{equation} 
In view of the ``fibrewise'' description of $\widehat{G}$ it is useful 
to generalise the wPW notion to multipliers:
\begin{defn}
\label{def:wpw}
 Let $G$ be a locally compact group and $\omega$ a type I multiplier
 on $G$. Then $G$ has $\omega$-wPW, if for every nonzero $f \in
 {\rm L}_c^\infty(G)$, ${\mathcal F}_{\omega}(f)$ does not vanish on a 
 set of positive $\omega$-Plancherel measure. 
\end{defn}
Now we have collected enough terminology to formulate the main result
of this paper. 
\begin{thm}
\label{thm:main}
 Let $G$ be type I, $N\lhd G$ a type I regularly embedded normal subgroup,
 with the additional property that for almost $\gamma \in \widehat{N}$, 
 $G/G_\gamma$ carries an invariant measure.
\begin{enumerate}
\item[(a)]
 Assume that, for almost all $\gamma \in \widehat{G}$,
 $G_\gamma /N$ has $\overline{\omega_\gamma}$-wPW. Moreover
 assume that the following condition holds:
 \begin{equation} \label{eqn:inv_wPW} \left. \begin{array}{c}
 \forall \varphi \in {\rm L}_c^\infty(G)  \\ ~\forall \mbox{ G - invariant }
 \Gamma \subset \widehat{N} \end{array} \right\}
 \left( \widehat{f}|_{\Gamma} = 0 \mbox{ and }
 \nu_N(\Gamma) > 0 \Rightarrow f=0 \right)
 ~~.\end{equation}
 Then $G$ has wPW. 
\item[(b)] Conversely, if $G$ has wPW, then (\ref{eqn:inv_wPW}) holds. 
\end{enumerate}
\end{thm}

The existence of invariant measures on $G/G_\gamma$ is ensured if
$G/N$ is abelian or compact. 
Note that relation $(\ref{eqn:inv_wPW})$
holds in particular when $N$ has wPW. 

\section{Proof of Theorem \ref{thm:main} }

The proof uses ideas and techniques very similar to 
those in \cite{KlLi}.
It rests on the interplay of several measure decompositions:
Of Haar measure on $G$ along shifts of subgroups on the group side, 
and of the Plancherel measures $\nu_G$ and $\nu_N$ according to 
(\ref{eqn:Pl_meas}) and (\ref{eqn:meas_decomp_Pl}).

For explicit calculations it is convenient to realize induced representations
${\rm Ind}_H^G \sigma$ on ${\rm L}^2(G/H; {\mathcal H}_\sigma)$.
Then the integrated representation acts on this space via operator-valued 
kernels. The proof of Theorem \ref{thm:main} uses explicit
formulas for these kernels, and their relationship to Fourier transforms 
of restrictions of $\varphi$ to cosets mod $N$. 
In the following, the restriction of a map
$f$ to a subset $Y$ of its domain is denoted by $f|_Y$. 
$f|_Y = 0$ is to be understood in the sense of vanishing
almost everywhere (with respect to a measure that is clear from the
context).

First let us recall a 
few basic results concerning cross-sections,
quasi-invariant measures and measure decompositions. If $H< G$, then
a cross-section $\alpha : G/H \to G$ is a measurable mapping fulfilling
$\alpha(\xi) H = \xi$, for all $\xi = g H \in G/H$. 
A cross-section $G/H \to G$ is called {\bf regular} if images of
compact subsets are relatively compact. All cross-sections in this
paper are assumed to be regular, which is justified by the following
lemma.
\begin{lemma}
 If $G$ is second countable and $H < G$ closed, there exists
 a regular cross-section.
\end{lemma}
\begin{prf}
Denote by $q : G \to G/H$ the quotient map. By 
\cite[Lemma 1.1]{Ma} there exists a Borel set $C$ of representatives
mod $H$, such that in addition for all $K \subset G$ compact
the set $C \cap q^{-1}(q(K))$ is relatively compact.

The associated cross-section is constructed by observing that 
$q|_C : C \to G/H$ is bijective. Then $q|_C$ is a measurable bijection 
between standard Borel spaces, hence the inverse map $\alpha$ is
also Borel, and it is a cross-section. Moreover, given any compact
$K \subset G/H$, there exists a compact subset $K_0 \subset K$
with $q(K_0)= K$ \cite[Lemma 2.46]{Fo}. 
Then $\alpha(K) = C \cap q^{-1}(q(K_0))$ is relatively compact.
\end{prf}

Given a cross-section, we may parametrise $G$ by the map 
\begin{equation} \label{eqn:par_G_1}
 H \times G/H \ni (h,\xi) \mapsto h \alpha(\xi)^{-1} ~~.
\end{equation}
This particular choice of parametrisation seems a bit peculiar, since
it refers to right cosets rather than left ones. Its benefit
will become apparent in the proof of Lemma \ref{lem:induced}.

Let us first take a closer look at the form that left Haar measure on $G$
takes in the parametrisation  (\ref{eqn:par_G_1}). We assume that $G/H$
carries an invariant measure, denoted in the following by $d\xi$.

\begin{lemma} \label{lem:meas_decomp}
 For all $f \in {\rm L}_c^\infty(G)$, 
\[
 \int_G f(x) dx = \int_{G/H} \int_H f(h \alpha(\xi)^{-1})
 \Delta_G(\alpha(\xi))^{-1} dh d\xi~~.
\]
\end{lemma}

\begin{prf}
\begin{eqnarray*} 
 \int_G f(x) dx & = & \int_G f(x^{-1}) \Delta_G(x^{-1}) dx \\
 & = & \int_{G/H} \int_H f(h^{-1} x^{-1}) \Delta_G(h^{-1} x^{-1})  dh
 d(xH) \\
 & = & \int_{G/H} \int_H f(h^{-1} \alpha(\xi)^{-1}) \Delta_G(h^{-1} 
 \alpha(\xi)^{-1}) dx d\xi \\
 & = & \int_{G/H} \int_H f(h \alpha(\xi)^{-1}) 
 \Delta_G(h \alpha(\xi)^{-1}) \Delta_H(h^{-1}) dh d\xi \\
 & = & \int_{G/H} \int_H f(h \alpha(\xi)^{-1}) \Delta_G(\alpha(\xi))^{-1}
 dh d\xi ~~,
\end{eqnarray*}
where we used Weil's integral formula and $\Delta_H = \Delta_G|_H$.
\end{prf}

Next we note a few technical details concerning the behaviour 
of restrictions of an ${\rm L}_c^\infty$-function to shifts of subgroups.
For this purpose one further piece of notation is necessary: 
The left- and right translation operators on $G$, denoted by
$R_x, L_y$, act via
\[ (L_y R_x f) (g) = f(ygx) ~~.\]

\begin{lemma} \label{lem:L1_norm_restr}
 Let $\varphi \in {\rm L}_c^\infty(G)$ and $H<G$; let $\alpha: G/H \to
G$ be a cross-section mapping compact sets to relatively 
compact sets. Consider the mapping
\[
 C : G/H \times G/H \to \RR ~~,~~ (\xi,\xi') \mapsto \| \left( L_{\alpha(\xi)}
 R_{\alpha(\xi')^{-1}} \varphi \right) |_H \|_1 
\]
where $\| \cdot \|_1$ denotes the ${\rm L}^1$-norm on $H$. 
\begin{enumerate}
\item[(a)] Given a compact set $K \subset G/H$, the set 
\[ \{ \xi' \in G/H : C(\xi,\xi') \not= 0 ~~\mbox{ for some }
 \xi \in K \} \]
is relatively compact.
\item[(b)] $C$ is bounded on compact subsets of $G/H \times G/H$. 
\end{enumerate}
\end{lemma} 
\begin{prf}
 For part $(a)$ note that $ \left( L_{\alpha(\xi)}
 R_{\alpha(\xi')^{-1}} \varphi \right) |_H$ is not
identically zero iff \\ $H \cap \alpha(\xi)^{-1} ({\rm supp}(\varphi)) \alpha
(\xi') \not= \emptyset$. Solving for $\alpha(\xi')$ we obtain the 
necessary condition 
\[ \alpha(\xi') \in ({\rm supp}(\varphi))^{-1} \alpha(K) H \]
or, equivalently
\[ \xi' \in  \left( ({\rm supp}(\varphi))^{-1} \alpha(K) H \right)/H~~.
\] 
By assumption on $\alpha$, $\alpha(K)$ is relatively compact, 
hence $(a)$ is shown.

For part $(b)$, it is enough to obtain an upper estimate for the
Haar measure of 
${\rm supp}( \left( L_{\alpha(\xi)}
 R_{\alpha(\xi')^{-1}} \varphi \right) |_H )$, 
and to consider compact sets of the form $K_1 \times K_2$. Here we
see that
\begin{eqnarray*}
 {\rm supp}( \left( L_{\alpha(\xi)}
 R_{\alpha(\xi')^{-1}} \varphi \right) |_H) ) & \subset &
 H \cap \alpha(\xi)^{-1} {\rm supp}(\varphi) \alpha(\xi') \\
 & \subset & \alpha(K_1^{-1})  {\rm supp}(\varphi) \alpha(K_2) \end{eqnarray*}
and the latter set is relatively compact.
\end{prf}

Now we can compute the action of the induced representation. 
\begin{lemma}
 \label{lem:induced}
 Let $G$ be a locally compact group, $H< G$ closed and $\sigma$
 a representation of $H$ acting on a separable Hilbert space
 ${\mathcal H}_\sigma$. 
 Let $\alpha: G/H \to G$ be a Borel cross-section and assume that
 there exists an invariant measure on $G/H$. We 
 realise $\pi = {\rm Ind}_H^G \sigma$ on the corresponding
 vector-valued space ${\rm L}^2(G/H;
 {\mathcal H}_\sigma)$ using $\alpha$. Then $\pi$ acts via
 \begin{equation}
 \label{eqn:ind_x}
 \pi(x) f(\xi) = \sigma\left( \alpha(\xi)^{-1} x \alpha(x^{-1} \xi) \right)
  f(x^{-1} \xi) ~~.
 \end{equation}
 For $\varphi \in {\rm L}_c^\infty(G)$,  $\pi(\varphi)$ acts 
 via
 \begin{equation} \label{eqn:op_kernel_1}
 \left[ \pi(\varphi) f \right] (\xi) = \int_{G/H} \Phi(\xi,\xi') 
 f(\xi') d\xi' \end{equation}
 where the right hand side converges in the weak sense for all $\xi \in
 G/H$, and $\Phi$ is an operator-valued integral kernel given by 
 \begin{equation} \label{eqn:op_kernel_2}
 \Phi(\xi,\xi') = \sigma\left( \left( L_{\alpha(\xi)}
 R_{\alpha(\xi')^{-1}} \varphi \right) |_H \right) \cdot
 \Delta_G(\alpha(\xi'))^{-1}
 ~~.
 \end{equation}
Moreover, we have the equivalence
\begin{equation} \label{eqn:ker_vanish}
 \pi(\varphi) = 0 \Leftrightarrow \Phi(\xi,\xi') = 0 ~~(a.e.)
\end{equation}
\end{lemma}

\begin{prf}
 Formula (\ref{eqn:ind_x}) is well-known.
 For weak convergence of the right-hand side of (\ref{eqn:op_kernel_1})
 let $\eta \in {\mathcal H}_\sigma$, and compute
\begin{eqnarray*}
 \lefteqn{\int_{G/H} \left|\langle \sigma\left( \left( L_{\alpha(\xi)}
 R_{\alpha(\xi')^{-1}} \varphi \right) |_H \right) 
 \Delta_G(\alpha(\xi'))^{-1} f(\xi'), \eta \rangle \right| d\xi' \le }\\
 & \le &  \int_{G/H} \left\| \sigma\left( \left( L_{\alpha(\xi)}
 R_{\alpha(\xi')^{-1}} \varphi \right) |_H \right) \right\|_\infty 
 \Delta_G(\alpha(\xi'))^{-1}  ~\| f(\xi') \| ~\|\eta \|~
 d\xi' \\
 & \le & \int_{G/H}  \left\|  \left( L_{\alpha(\xi)}
 R_{\alpha(\xi')^{-1}} \varphi \right) |_H
 \right\|_1  \Delta_G(\alpha(\xi'))^{-1}  ~\|f(\xi')\| d\xi' ~\| \eta \|~~. 
\end{eqnarray*}
By Lemma \ref{lem:L1_norm_restr} (a) the map
\[
 \xi' \mapsto  \left\|  \left( L_{\alpha(\xi)}
 R_{\alpha(\xi')^{-1}} \varphi \right) |_H
 \right\|_1    \Delta_G(\alpha(\xi'))^{-1} 
\]
is compactly supported, and also bounded, by Lemma \ref{lem:L1_norm_restr} (b)
and boundedness of $\Delta_G^{-1} \circ \alpha$ on the support. 
Hence the map is square-integrable, and
an application of the Cauchy-Schwarz inequality finishes the proof of
weak convergence. 

For the integrated transform, let $f,g \in {\rm L}^2(G/H;
{\mathcal H}_\sigma)$.  Then, by (\ref{eqn:ind_x}), the weak definition 
of $\pi(\varphi)$ yields
 \begin{eqnarray}
\lefteqn{ \langle \pi(\varphi) f , g \rangle = \nonumber
 }\\  & = & \int_{G/H}
 \int_G   \varphi(x) \left\langle \sigma\left( \alpha(\xi)^{-1} x
 \alpha(x^{-1} \xi) \right) f(x^{-1} \xi), g(\xi) \right\rangle
 dx d\xi \nonumber \\
 & = & \int_{G/H}
 \int_G \varphi(x) \left\langle \sigma\left( \alpha(\xi)^{-1} x
 \alpha(x^{-1} \xi) \right)  f(x^{-1} \alpha(\xi) H), g(\xi)
 \right\rangle dx d\xi \nonumber \\ 
 & = & \label{eqn:shift}  \int_{G/H}
 \int_G \varphi(\alpha(\xi)x)   
 \left\langle \sigma\left( x
 \alpha((\alpha(\xi)x)^{-1} \xi) \right)  f(x^{-1} H), g(\xi)
 \right\rangle dx d\xi \\
 & = & \nonumber
 \int_{G/H} \int_{G/H} \int_H \varphi(\alpha(\xi) h \alpha(\xi')^{-1})
 \Delta_G(\alpha(\xi'))^{-1} \\
 & &  \label{eqn:final} \hspace*{2cm} 
  \left\langle \sigma( h )  f(\xi'), g(\xi)
 \right\rangle dh d\xi' d\xi ~~.
\end{eqnarray}
 Here (\ref{eqn:shift}) was obtained by a left translation in the
 integration variable $x$. (\ref{eqn:final}) used
 Lemma \ref{lem:meas_decomp}, as well as the calculations
\[ 
 \alpha(\xi') h^{-1} H = \xi'
\]
and 
\begin{eqnarray*}
  h \alpha(\xi')^{-1} \alpha( (\alpha(\xi) h \alpha(\xi')^{-1} )^{-1} \xi ) & = & 
 h \alpha(\xi')^{-1} \alpha( 
 \alpha(\xi') h^{-1} \alpha(\xi)^{-1} \xi) \\ & = & h \alpha(\xi')^{-1}
 \alpha(\alpha(\xi') h^{-1} H) \\
 &  = & h \alpha(\xi')^{-1} \alpha(\alpha(\xi') H) \\
 & = & h ~~.\end{eqnarray*}
Using the definition of weak integrals, we may continue from 
(\ref{eqn:final}) to obtain 
 \begin{eqnarray*}
\lefteqn{ \langle \pi(\varphi) f , g \rangle = \nonumber
 }\\ & = &
 \int_{G/H} \int_{G/H}
 \langle \Delta_G(\alpha(\xi'))^{-1}  \sigma\left( \left( L_{\alpha(\xi)}
 R_{\alpha(\xi')^{-1}} \varphi \right)|_H \right)
  f(\xi') , g(\xi) \rangle d\xi' d\xi \\
 & = & \int_{G/H} \left\langle \int_{G/H} \Phi(\xi,\xi') f(\xi') d\xi', g(\xi) \right\rangle
 d\xi ~~,
\end{eqnarray*}
hence the pointwise definition $(\ref{eqn:op_kernel_1})$ indeed coincides
with $\pi(\varphi)$. 

Now the direction ``$\Leftarrow$'' of (\ref{eqn:ker_vanish}) is immediate.
For the other direction assume that $\Phi$ does not vanish almost
everywhere. Pick an ONB $(\eta_i)_{i \in I}$ of ${\mathcal H}_\sigma$;
since ${\mathcal H}_\sigma$ is separable, $I$ is countable. 
Since 
\[ 
 \Phi(\xi, \xi') = 0 \Leftrightarrow \forall i,j \in I : \langle \Phi(\xi,\xi')
 \eta_i, \eta_j \rangle  = 0
\]
there exists a pair $(i,j)$ and a set $A \subset G/H \times G/H$ of 
positive measure such that 
\[
  \langle \Phi(\xi,\xi') \eta_i, \eta_j \rangle \not= 0 
\]
for all $(\xi,\xi') \in A$. By passing to a smaller set we may 
assume that in addition
$A \subset G/H \times K$ for a compact $K \subset G/H$ (observing
that $G/H$ is $\sigma$-compact). 

Now define the auxiliary operator $T: {\rm L}^2(X) \to {\rm L}^2(X)$ by 
\[
 (Tf)(\xi) = \chi_K (\xi) \langle \sigma(\varphi) (f \cdot \eta_i)(\xi), \eta_j
 \rangle~~.
\] 
 Then $T$ is an integral operator with kernel
\[
 (\xi,\xi') \mapsto \chi_K(\xi) 
 \langle \Phi(\xi,\xi') \eta_i, \eta_j \rangle~~,
\]
which by construction is nonzero.
By Lemma \ref{lem:L1_norm_restr} this kernel is bounded and compactly
supported, hence in ${\rm L}^2(G/H\times G/H)$. But for this space the
map from kernel to integral operator is a unitary operator onto the 
space of Hilbert-Schmidt operators on ${\rm L}^2(G/H)$; in particular 
the map is one-to-one.
Hence $T\not=0$, which implies $\sigma(\varphi) \not= 0$. 
\end{prf}

\noindent{\bf Proof of Theorem \ref{thm:main}.}
 Let $\varphi \in {\rm L}_c^\infty(G)$ be given with
 $\pi(\varphi) = 0$ for $\pi$ in a set of positive Plancherel measure. 
 Then, by (\ref{eqn:Pl_meas}), there exists a $G$-invariant subset
 $\Gamma \subset \widehat{N}$ of positive Plancherel measure and
 subsets $B_{G. \gamma} \subset \widehat{G_\gamma/N}^{\overline{\omega_\gamma}}$
 ($\gamma \in \Gamma$)  with
 \[ \nu_{G_\gamma/N,{\overline{\omega_\gamma}}}(B_{G. \gamma})> 0 \mbox{ and }
 \pi_{G.\gamma, \sigma}(\varphi) = 0 \mbox{, 
 for all } \sigma \in B_{G. \gamma}.\]

 Now fix $\gamma \in \Gamma$. Our aim is to relate the equation 
 $\pi_{G.\gamma,\sigma}(\varphi)=0$ for 
 $\sigma$ in a set of positive projective Plancherel measure in
 $\widehat{G_\gamma/N}^{\overline{\omega_\gamma}}$ to certain 
 Fourier transforms, using the integral kernel calculus.
 For this purpose we use
 Borel cross-sections $\alpha: G/G_{\gamma} \to G$ and
 $\vartheta : G_{\gamma} / N \to G_\gamma$.
 In the following calculations, all quotients carry invariant measures.
 We also need the continuous homomorphism $\delta: G \to (\RR^+,\cdot)$ 
 defined by 
 picking $B \subset N$ of positive finite measure and letting
 $\delta(x) = \frac{|B|}{|x B x^{-1}|}$.

 By Lemma \ref{lem:induced},  $\pi_{G.\gamma,\sigma}(\varphi)$
 has the operator-valued kernel
 \begin{eqnarray*}
 \lefteqn{ \Phi(\xi,\xi') = } \\ & = & \int_{G_{\gamma}} \left( \gamma' (y) \otimes \sigma(y)
 \right) \varphi(\alpha(\xi) y \alpha(\xi')^{-1}) dy ~ \Delta_G(\alpha(\xi'))^{-1} \\
 & = & \int_{G_\gamma/N} \int_N \left( \gamma' \left(n \vartheta(h)^{-1}
 \right) \otimes \sigma (h)^{-1} \right)  \varphi(\alpha(\xi) n \vartheta(h)^{-1}
 \alpha(\xi')^{-1}) \\
 & & \hspace*{3cm}
 \Delta_{G_{\gamma}}(\vartheta(h))^{-1}
 dn dh ~ \Delta_G(\alpha(\xi'))^{-1} \\
 & = & \int_{G_\gamma/N} \left( \int_N \gamma(n) \varphi(\alpha(\xi) n
 \vartheta (h)^{-1} \alpha(\xi')^{-1}) dn \circ \gamma'(\vartheta(h)^{-1}) 
 \Delta_{G_\gamma}(\vartheta(h))^{-1} \right) \\
 & & \hspace*{2cm} \otimes \sigma(h^{-1}) dh ~ \Delta_G(\alpha(\xi'))^{-1} \\
 & = & \int_{G_\gamma/N} F_{\xi,\xi'}(h) \otimes \sigma(h^{-1})
 dh  ~\Delta_G(\alpha(\xi'))^{-1} ~~,
 \end{eqnarray*}
where 
\begin{eqnarray} \nonumber 
\lefteqn{ F_{\xi,\xi'}(h) =} \\
 &=& \nonumber \int_N \gamma(n) \varphi(\alpha(\xi) n
 \vartheta (h)^{-1} \alpha(\xi')^{-1}) dn \circ \gamma'(\vartheta(h)^{-1})
  \Delta_{G_\gamma}(\vartheta(h))^{-1} \\
 & = & \nonumber \delta(\alpha(\xi)) 
 \int_N \gamma\left(\alpha(\xi)^{-1} n \alpha(\xi) \right) ~
 \varphi(n \alpha(\xi) \vartheta (h)^{-1} \alpha(\xi')^{-1}) dn
 \circ \gamma'(\vartheta(h)^{-1})  \Delta_{G_\gamma}(\vartheta(h))^{-1}\\
\label{eqn:ker_explicit}
 & = & \delta(\alpha(\xi)) \left[
 \left( \alpha(\xi). \gamma\right) \left( \left( R_{\alpha(\xi) 
 \vartheta(h)^{-1} \alpha(\xi')^{-1}} \varphi \right)|_{N} \right) \right]
 \circ \gamma'(\vartheta(h)^{-1}) \Delta_{G_\gamma}(\vartheta(h))^{-1} ~~.
\end{eqnarray}
Here it is important to note that for fixed $(\xi,\xi')$,
the operator-valued function $F_{\xi,\xi'}$ 
has compact support: A short calculation establishes that 
\[ \left( R_{\alpha(\xi) 
 \vartheta(h)^{-1} \alpha(\xi')^{-1}} \varphi \right)|_{N} \not=0 \]
only if $h \in \alpha(\xi')^{-1} ({\rm supp} (\varphi))^{-1} \alpha(\xi)^{-1}
N =: K_0$, and $K_0$ is a compact subset of
$G/N \supset G_\gamma/N$. Moreover, for $h \in K_0$
\begin{eqnarray*}
\lefteqn{ \left\| \delta(\alpha(\xi)) \left( \alpha(\xi). \gamma\right)
\left( \left( R_{\alpha(\xi) 
 \vartheta(h)^{-1} \alpha(\xi')^{-1}} \varphi \right)|_{N} \right)
 \circ \gamma'(\vartheta(h)^{-1})  \Delta_{G_\gamma}(\vartheta(h))^{-1}
 \right\|_\infty}\\ & \le &\delta(\alpha(\xi)) 
 \left\| \left( R_{\alpha(\xi) 
 \vartheta(h)^{-1} \alpha(\xi')^{-1}} \varphi \right)|_{N} \right\|_1
  \Delta_{G_\gamma}(\vartheta(h))^{-1} \\
  & \le & \delta(\alpha(\xi)) \| \varphi
 \|_\infty \left| N \cap {\rm supp}(\varphi) \alpha(\xi')
 \vartheta(K_0) \alpha(\xi)^{-1} \right|  \Delta_{G_\gamma}(\vartheta(h))^{-1}
 ~~.
\end{eqnarray*}
The middle term is the measure of a fixed relatively compact subset of $N$, 
by regularity of $\vartheta$, and
the last term is bounded on the compact support.
Hence the map $h \mapsto \| F_{\xi,
\xi'} (h) \|_\infty$ is in ${\rm L}_c^\infty(G_\gamma/N)$. 

Now, for fixed $\gamma \in \Gamma$ and $\sigma \in B_{G.\gamma}$,
relation (\ref{eqn:ker_vanish}) and $\Delta_G > 0$ imply that 
\begin{equation} \label{eqn:kernel_vanish_1}
 \int_{G_\gamma/N} F_{\xi,\xi'}(h) \otimes \sigma(h^{-1}) dh = 0 
\end{equation}
for almost every $(\xi,\xi')$, where the set of these $(\xi,\xi')$
may depend on $\sigma$. However, an application of Fubini's Theorem
provides a conull subset of $C \subset G/G_\gamma \times G/G_\gamma$,
such that (\ref{eqn:kernel_vanish_1}) holds for  all $(\xi,\xi') \in C$
and all $\sigma$ in a conull subset of $B_{G.\gamma}$, possibly depending
on $(\xi,\xi')$. Now fix $(\xi,\xi') \in C$,
as well as ONB's $(\eta_i)_{i \in I} \subset {\mathcal H}_\gamma$,
 $(\beta_j)_{j \in J} \subset {\mathcal H}_\sigma$.
It follows that 
\begin{eqnarray*}
 0 & = & \langle \Phi_{\xi,\xi'}
 (\eta_i \otimes  \beta_j), (\eta_k \otimes \beta_\ell) \rangle  \\ 
& = &  \int_{G_\gamma/N} \langle F_{\xi,\xi'} (h) \eta_i, \eta_k \rangle
 \langle \beta_j, \sigma(h)^{-1} \beta_\ell
 \rangle dh \\
& = & \langle \sigma(\Psi_{i,k}) \beta_j, \beta_\ell \rangle~~,
\end{eqnarray*}
where we used that $dh$ is left Haar measure on $G_\gamma/N$, and the
scalar-valued function $\Psi_{i,k}$ given by
\[ \Psi_{i,k}(h) =  \langle \eta_i, F_{\xi,\xi'} (h) \eta_k \rangle~~,~~
\mbox{ satisfying } |\Psi_{i,k}(h)| \le
\| F_{\xi,\xi'} (h) \|_\infty \|\eta\|~\|\eta'\|~~.\]
In particular $\Psi_{i,k} \in {\rm L}_c^\infty(G_\gamma/N)$.
Thus we are finally in a position to use the assumption that
$G_\gamma /N$ has $\overline{\omega_\gamma}$-wPW, yielding for all $i,k$
that
$\Psi_{i,k}=0$ on a joint conull subset. But this clearly entails 
$F_{\xi,\xi'}(h) = 0$ almost everywhere.

Since obviously 
\[ F_{\xi,\xi'} (h) = 0 \Leftrightarrow 
 \left( \alpha(\xi). \gamma\right) \left( \left( R_{\alpha(\xi)
 \vartheta(h)^{-1} \alpha(\xi')^{-1}} \varphi \right)|_N \right) = 0
\]
our considerations so far have established for fixed $\gamma
\in \Gamma$, that for almost all $(\xi, \xi',h) \in G/G_\gamma \times 
 G/G_\gamma \times G_\gamma /N$ 
\begin{equation} \label{eqn:ker_vanish_ae}
  \left( \alpha(\xi). \gamma\right) \left( \left( R_{\alpha(\xi)
 \vartheta(h)^{-1} \alpha(\xi')^{-1}} \varphi \right)|_N \right) = 0
\end{equation}

Now choose a measurable cross-section $\Lambda : G/N \to G$. By definition 
of the Borel structure on $\widehat{N}$, the map 
\[
 \widehat{N} \times G/N \ni (\gamma, s) \mapsto \left\| \gamma
 \left( \left(R_{\Lambda(s)^{-1}} \varphi \right)|_N \right) \right\|_\infty   
\]
is Borel. 
Moreover, for a fixed $\xi \in G/G_\gamma$, the set 
 \[ \{ \alpha(\xi) \vartheta(h)^{-1} \alpha(\xi')^{-1} :
 h \in G_\gamma/N, \xi' \in  G/G_\gamma \} \]
 is a set of representatives of $G/N$, since $N$ is normal.
 In particular, for every $s \in G/N$ there exists $(\xi',h) \in
 G/G_\gamma \times G_\gamma /N$ such that $N \Lambda(s)^{-1}
 = N \alpha(\xi) \vartheta(h)^{-1} \alpha(\xi')^{-1}$.
 Hence $\Lambda(s)^{-1} = n \alpha(\xi) \vartheta(h)^{-1} \alpha(\xi')^{-1}$,
 for suitable $n \in N$. Thus (\ref{eqn:ker_vanish_ae}) implies that 
 for fixed $\gamma \in \Gamma$ and almost all $\xi \in G/G_\gamma$ the set
\[ \{ s \in G/N :  \left\| (\alpha(\xi).\gamma)
 \left( \left(R_{\Lambda(s)^{-1}} \varphi \right)|_N \right) \right\|_\infty
  = 0 \} \]
 has a complement of measure zero.

 On the other hand, $\xi \mapsto \alpha(\xi).\gamma$ yields a
 bijection between $G/G_\gamma$ and $G.\gamma$, and the image
 of the invariant measure on $G/G_\gamma$ is 
 equivalent to the measure $\mu_{G.\gamma}$ 
 appearing in (\ref{eqn:meas_decomp_Pl}).
 Summarising, we obtain that
\begin{eqnarray*}
0 & = &  \int_{\Gamma/N} \int_{G . \gamma} \int_{G/N}
 \|  \gamma \left( \left(R_{\Lambda(s)^{-1}} \varphi \right)|_N \right) 
 \|_{\infty} ~ds~d\mu_{G.\gamma} (\gamma)~d\overline{\nu}(G.\gamma) \\
 & = & \int_{G/N} \int_\Gamma  \|  \gamma
  \left( \left(R_{\Lambda(s)^{-1}} \varphi \right)|_N \right)
 \|_{\infty} d\nu_N
 (\gamma) ds ~~.
\end{eqnarray*}
Here the second equation uses the measure decomposition
(\ref{eqn:meas_decomp_Pl}) and Fubini's Theorem.
But the latter integral implies for almost all $s \in G/N$, that
$  \left( \left(R_{\Lambda(s)^{-1}} \varphi \right)|_N \right)^\wedge = 0$ 
on the $G$-invariant subset $\Gamma \subset \widehat{N}$. 
Since $ \left(R_{\Lambda(s)^{-1}} \varphi \right)|_N \in {\rm L}_c^\infty(N)$,
an appeal to assumption 
(\ref{eqn:inv_wPW}) yields $ \left(R_{\Lambda(s)^{-1}} \varphi \right)|_N =0$
for almost every coset $s$.
Hence $\varphi=0$, which finishes the proof of $(a)$.

For the proof of $(b)$, assume that $\varphi_0 \in {\rm L}_c^\infty(N) 
\setminus \{ 0 \}$
is a counterexample to (\ref{eqn:inv_wPW}), i.e. there exists
a $G$-invariant $\widetilde{\Gamma} \subset \widehat{N}$ of positive measure 
such that $\widehat{\varphi_0}$ vanishes on $\Gamma$. 
Let $K \subset G/N$ be some compact set of positive measure, 
and let $\Lambda: G/N \to G$ be a cross-section.
Then
\[
 \varphi(n \Lambda(s)^{-1}) = \varphi_0(n) \chi_K(s) 
\]
defines a nonzero $\varphi \in {\rm L}_c^\infty(G)$. Let $\Sigma
= \{ \pi \in \widehat{G} : \widehat{\varphi}(\pi) = 0 \}$. 
We intend to show $\pi_{G.\gamma,\sigma} \in \Sigma$ for all
$\gamma \in \widetilde{\Gamma}$
and all $\sigma \in \widehat{G.\gamma}^{\overline{\omega_\gamma}}$.
By equation (\ref{eqn:ker_explicit}) this amounts to proving
 \[
 0 = \left( \alpha(\xi). \gamma\right)
 \left( \left( R_{\alpha(\xi)
 \vartheta(h)^{-1} \alpha(\xi')^{-1}} \varphi \right)|_N \right)~~,
\]
where $\alpha, \vartheta$ are cross-sections associated to $G_\gamma/N$
and $G/G_\gamma$.
Now for $n \in N$,
\begin{eqnarray*}
\left( R_{\alpha(\xi)
 \vartheta(h)^{-1} \alpha(\xi')^{-1}} \varphi \right)|_N (n)
 & = & \varphi(n  \alpha(\xi) \vartheta(h)^{-1} \alpha(\xi')^{-1}) \\
 & = & \varphi(n n' \Lambda(s)^{-1}) \\
 & = & \varphi_0(n n') \chi_K(s)
\end{eqnarray*}
where $s = \alpha(\xi') \vartheta(h) \alpha(\xi)^{-1} N$ and
$n' \in N$ is suitably chosen, and independent of $n$.
Since $\alpha(\xi).\gamma \in \Sigma$, it follows that 
\[
 \left( \alpha(\xi). \gamma\right)
 \left( \left( R_{\alpha(\xi)
 \vartheta(h)^{-1} \alpha(\xi')^{-1}} \varphi \right)|_N \right)
  = \chi_K(s) \left( \alpha(\xi). \gamma\right) (\varphi_0) \circ
  \left( \alpha(\xi). \gamma\right) (n') = 0
 \]
Hence we can compute 
 \begin{eqnarray*}
 \nu_G(\Sigma) & = & \int_{\widehat{N}/G} \int_{\widehat{G_\gamma/N}^{\overline{\omega_\gamma}}} 
 \chi_{\Sigma} (\pi_{G.\gamma, \sigma}) ~
 d\nu_{G_\gamma/N,\overline{\omega_\gamma}} (\sigma)~ d\overline{\nu}
 (G.\gamma) \\
 & \ge &  \int_{\widetilde{\Gamma}/G} \nu_{G_\gamma/N,\overline{\omega_\gamma}}
 (\widehat{G_\gamma/N}^{\overline{\omega_\gamma}})~
  d\overline{\nu}  (G.\gamma) \\
 & > & 0~~,
 \end{eqnarray*}
therefore $\varphi$ is the desired counterexample to wPW on $G$.
 {\hfill $\Box$}

\section{Applications and Examples}

In this section we apply Theorem \ref{thm:main} to a variety of 
cases, and discuss the necessity of its assumptions.
Unless otherwise stated, our standing assumptions are: $G$ is
second countable, 
$G$ and $N \lhd G$ are of type I,
with $N$ regularly embedded.
\begin{cor}
\label{cor:free}
 Assume that the dual action of $G/N$ is free $\nu_N$-almost everywhere.
 Then $G$ has wPW iff condition $(\ref{eqn:inv_wPW})$ holds. 
\end{cor}

\begin{cor} \label{cor:nonunimod}
 Let $G$ be nonunimodular, and $N = {\rm Ker}(\Delta_G)$. Then $G$
 has wPW iff $(\ref{eqn:inv_wPW})$ holds; in particular, if $N$ has
 wPW. 
\end{cor}
\begin{prf}
 This is an immediate consequence of Corollary \ref{cor:free},
 which applies to this setting by \cite[Section 5]{DuMo}.
\end{prf}

For split compact extensions the freeness of the operation turns out
to be necessary also. Since a compact group has wPW iff it is trivial,
the next theorem provides a class of extensions for which the conditions 
of Theorem \ref{thm:main} are necessary and sufficient. 
\begin{thm} \label{thm:motion}
 Assume that $G = N \semdir K$, with $K$ compact. 
 Then $G$ has wPW iff condition (\ref{eqn:inv_wPW}) holds
 and the dual action of $G/N$ is free $\nu_N$-almost everywhere. 
\end{thm}
\begin{prf}
 The ``if''-part is Corollary \ref{cor:free}. For the ``only-if''-part,
 necessity of (\ref{eqn:inv_wPW}) was noted in Theorem \ref{thm:main}.
 We denote elements of $G$ by pairs $(n,k) \in N \times K$, and the conjugation
 action of $K$ on $N$ by $K \times N \ni (k,n) \mapsto k.n$.
 Define the little fixed groups as $K_\gamma = G_\gamma \cap K$;
 since $G$ is a semidirect product, $G_\gamma = N \semdir K_\gamma$.

 We define
 \[
 \widetilde{\Sigma} = \{ \gamma \in \widehat{N} : K_\gamma \not= \{ 1 \} \}~~.
 \]
 Let us first show that $\widetilde{\Sigma}$ is a Borel subset of 
$\widehat{N}$. For this purpose consider the space $X$ of closed subgroups
 of $K$, endowed with the compact open
 topology. By \cite[Proposition II.2.3]{AuMo}, the stabilizer map $\widehat{N}
 \to X$ is Borel. Moreover, the complement of $\widetilde{\Sigma}$ is nothing
 but the inverse image of the trivial subgroup under the stabiliser
 map, hence Borel.
 Thus $\widetilde{\Sigma}$ is Borel also. Assuming that
 $\nu_N(\widetilde{\Sigma})>0$, we need to construct a $\varphi \in
 {\rm L}_c^\infty(G)$ such that $\widehat{\varphi}$ vanishes on a set of
 positive measure.

 Let $\varphi_0 \in {\rm L}_c^\infty (N) \setminus \{ 0 \}$ with
 $\varphi_0 \ge 0$ be given, and let 
 \[
 \varphi(n,k) = \varphi_1(n) = \int_K \varphi_0 (k'.n) dk'  
 \]
 which yields a nonzero element $\varphi \in {\rm L}_c^\infty(G)$.
 Next we show 
 \begin{equation} \label{eqn:four_vanish}
 \forall \gamma \in \widetilde{\Sigma}~,~ \forall \sigma \in 
 \widehat{K_\gamma} \setminus \{ 1_{K_\gamma} \} ~:~
 \pi_{G.\gamma,\sigma} (\varphi) = 0 \end{equation}
 where $1_{K_\gamma}$ denotes the trivial representation of $K_\gamma$.

 Since $G$ is a semidirect product, all Mackey obstructions are trivial.
 In addition, we can assume that all cross-sections arising below in fact map 
 into $K < G$, i.e. $\vartheta(h) = 
 (e_N, \widetilde{\vartheta}(h))$ etc.
 Moreover, since $K$ is compact all involved
 measures can be chosen invariant. In this setting the calculations
 from the proof of Theorem \ref{thm:main} yield that
 $\pi_{G.\gamma,\sigma}(\varphi)$ acts on
 ${\rm L}^2(K/K_\gamma;{\mathcal H}_\gamma \otimes {\mathcal H}_\sigma)$ via 
 \[
 \Phi(\xi,\xi') = \int_{K_\gamma} \left( \alpha(\xi).\gamma \right) 
 \left( \left( R_{\alpha(\xi) \vartheta(h)^{-1} \alpha(\xi')^{-1}} \varphi
 \right)|_N \right) \otimes 
 \sigma(h^{-1}) dh ~~.
 \]
 In order to prove that $\Phi$ vanishes, it is enough to show for all
  $\xi,\xi'$ that the map
 \[
 F_{\xi,\xi'} : h \mapsto  \left( \alpha(\xi).\gamma \right) 
 \left( \left( R_{\alpha(\xi) \vartheta(h)^{-1} \alpha(\xi')^{-1}} \varphi
 \right)|_N \right)  ~~. 
 \]
 is constant on $K_\gamma$. Note first that by construction of $\varphi$,
 and the fact that $\vartheta, \alpha$ map into $K$ that 
 \[  \left( \left( R_{\alpha(\xi) \vartheta(h)^{-1} \alpha(\xi')^{-1}} \varphi
 \right)|_N \right) (n) = 
 \varphi_1(n)
 \]
 is independent of $\xi,h, \xi'$ and invariant under the action of $K$. Hence we obtain
 \[ F_{\xi,\xi'} (h ) =  \left( \alpha(\xi).\gamma \right) 
 \left( \varphi_1 \right) = \gamma \left( \varphi_1 \right) ~~.\]
 Hence $F_{\xi,\xi'}$ is constant, and thus $\pi_{G.\gamma,\sigma}(\varphi) = 0$.

Hence, defining the Borel subset 
 \[
 \Sigma = \{ \pi \in \widehat{G} : \widehat{\varphi} (\pi) \not= 0 \}~~,
 \]
we can use (\ref{eqn:Pl_meas}) and (\ref{eqn:four_vanish}) to estimate 
 \begin{eqnarray*}
 \nu_G(\Sigma) & = & \int_{\widehat{N}/G} \int_{\widehat{K_\gamma}}
 \chi_{\Sigma} (\pi_{G.\gamma, \sigma}) d\nu_{K_\gamma} (\sigma) d\overline{\nu}
 (G.\gamma) \\
 & \ge &  \int_{\widetilde{\Sigma}/G} \int_{\widehat{K_\gamma}}
 \chi_{\Sigma} (\pi_{G.\gamma, \sigma}) d\nu_{K_\gamma} (\sigma) d\overline{\nu}
 (G.\gamma) \\ 
 & \ge &  \int_{\widetilde{\Sigma}/G} \nu_{K_\gamma} \left(\widehat{K_\gamma}
 \setminus \{ 1_{K_\gamma} \} \right) d\overline{\nu}(G.\gamma) \\
 & > & 0~~.
 \end{eqnarray*}
Here the last inequality is due to the fact that the integrand is strictly
positive on $\widetilde{\Sigma}/N$, and we assumed
$\nu_N(\widetilde{\Sigma})>0$. 
\end{prf}

Applying the theorem to motion groups yields that $G = \RR^n \semdir
SO(n)$ has wPW iff $n \le 2$. 
 
Another extreme case is given by an almost everywhere trivial action of
$G/N$ on $\widehat{G}$. The following corollary also covers direct
product groups.
\begin{cor}
 Assume that $G/N$ acts trivially $\nu_N$-almost everywhere. 
 If $G$ has wPW, then $N$ has wPW. Conversely, if both $N$ and $G/N$ have wPW,
 then so does $G$. 
\end{cor}
Note that $G/N$ need not have wPW, even if $G$ does: Simply take
$G=\RR$ and $N=\ZZ$.

A result similar to the following is formulated for the 
so-called {\em topological
Paley-Wiener condition} in \cite[Theorem 2.2]{KaLaSc}.
\begin{cor} \label{cor:abelian_ext}
 Suppose that $G/N$ is abelian and compact-free. Assume in addition 
 that either almost all Mackey obstructions vanish, or that $G/N$
 is compactly generated. Then, if condition \ref{eqn:inv_wPW} holds,
 $G$ has wPW.
\end{cor}
\begin{prf}
 Let us first deal with the case of vanishing Mackey obstructions.
 Recall that $G/N$ is compact-free iff it has no nontrivial compact
 subgroups. For abelian groups, this is equivalent to wPW by
 \cite[Theorem 3.2]{KaLaSc}. Moreover, if $G/N$ is compact-free, so are 
 all its closed subgroups; in particular, the little fixed groups also
 have wPW. Hence Theorem \ref{thm:main} implies wPW for $G$.

 If $G/N$ is compact-free and compactly generated, the structure theorem 
 for LCA groups yields $G/N \cong \RR^k \times \ZZ^ \ell$, and the
 little fixed groups have a similar structure. Hence 
 Theorem \ref{thm:main} together with the next lemma yield that
 $G$ has wPW.
\end{prf}

\begin{lemma}
 Let $G = \RR^k \times \ZZ^\ell$, and $\omega$ a type I multiplier 
 on $G$. Then $G$ has $\omega$-wPW.  
\end{lemma}

\begin{prf}
 We use the description of $\widehat{G}^\omega$ given in \cite{BaKl}.
 We may assume that $\omega$ is normalised. Then the map
 \[
 h_\omega : G \to \widehat{G}~~,~~ h_\omega (x) (y) = \omega(x,y)
 \overline{\omega(y,x)} 
 \]
defines a continuous homomorphism. Denote the kernel of this homomorphism
by $S_\omega$.  $\omega$ is called {\em totally skew}
if $S_\omega$ is trivial.
 By \cite[Theorem 3.1]{BaKl}, we may then assume that
$\omega$ is lifted from a totally skew cocycle $\omega_1$ of $G/S_\omega$.
Moreover, $\widehat{G/S_\omega}^{\omega_1} = \{ \pi \}$, and the mapping 
$\widehat{G} \ni \gamma \mapsto \gamma \pi \in \widehat{G}^\omega$
is continuous and onto (also by \cite[Theorem 3.1]{BaKl}). 
This map is constant on $S_\omega^\bot$, giving rise to 
a homeomorphism between $\widehat{G}^\omega$ and $\widehat{G}/S_\omega^\bot
\simeq \widehat{S_\omega}$.
Moreover, the projective Plancherel measure can be chosen as the 
Haar measure on $\widehat{G}/S_\omega^\bot$. To see this consider
the unitary action of $\widehat{G}$ on ${\rm L}^2(G)$ defined by pointwise
multiplication, $(M_\gamma f)(x) = \gamma(x) f(x)$. Then it is straightforward
to compute that on the (projective) Fourier transform side, $\widehat{G}$
acts via shifts, 
\[(\gamma_1 \pi) (M_{\gamma_2}f) = (\gamma_1 \gamma_2 \pi)
(f) ~~.\]

On the other hand, since $G$ is unimodular, there exists a choice
of $\nu_{G,\omega}$ such that ${\mathcal F}_\omega$ extends
to a unitary equivalence 
\[ {\rm L}^2(G) \to \int_{\widehat{G}^\omega}^\oplus 
{\rm HS}({\mathcal H}_\rho) d\nu_{G,\omega}(\rho) \simeq 
{\rm L}^2(\widehat{G}/S_\omega^\bot, d\nu_{G,\omega}) \otimes
{\rm HS}({\mathcal H}_\pi) ~~.\]
Since we already know that the shifts on 
$\widehat{G}/S_\omega^\bot$ yield unitary operators on 
${\rm L}^2(\widehat{G}/S_\omega^\bot,d\nu_{G,\omega})$, it follows
that $\nu_{G,\omega}$ is shiftinvariant.

Now assume that $f \in {\rm L}_c^\infty(G)$ is such that 
$\rho(f) = 0$ on a set of projective Plancherel measure zero.
Then the map $\widehat{G} \ni \gamma \mapsto (\gamma \pi) (f)$ 
also vanishes on a set of positive Plancherel measure. 
Pick an ONB $(\eta_i)_{i \in I}$ of ${\mathcal H}_\pi$.
Then for all $i,j \in I$
\begin{eqnarray*}
 0 & = & \langle (\gamma \pi)(f) \eta_i,\eta_j \rangle \\
 & = & \int_G \gamma(x) f(x) \langle \pi(x) \eta_i, \eta_j \rangle dx ~~,
\end{eqnarray*}
for all $\gamma$ in a set of positive measure. Hence wPW for $G$ implies that
\[ 0 = f(x) \langle \pi(x) \eta_i, \eta_j \rangle \]
for all $i,j \in I$ and almost all $x \in G$. On the other hand,
the fact that $\pi(x)$ is unitary implies for all $x \in G$ that
$\langle \pi(x) \eta_i, \eta_j \rangle \not= 0$ for some pair $(i,j)$.  
Thus we obtain $f = 0$ almost everywhere. 
\end{prf}

We will next show that 
iterated application of Corollary \ref{cor:abelian_ext}
allows to establish wPW for a large class of solvable Lie groups,
thus extending the results from \cite{LipRos,Gar1,ArnLud,Gar2}.
In the following, the term ``Lie group'' is shorthand for simply
connected, connected Lie group.
We first start with an observation that is probably folklore, and which
ensures that Corollary \ref{cor:abelian_ext} can be used iteratively.
We include a proof since we could not obtain a reference.
\begin{lemma}
\label{lemma:nilp_reg_emb}
 Let $G$ be an exponential Lie group and $N \lhd G$ a
 closed connected nilpotent normal subgroup. Then $N$ is regularly embedded. 
\end{lemma}
\begin{prf}
 Denote the Lie algebras of $G$,$N$ by
 ${\mathfrak g},{\mathfrak n}$, and by  ${\rm Ad}_G^*$ and ${\rm Ad}_N^*$ the
 coadjoint actions of $G$ and $N$ respectively. 

 $G$ being exponential implies that ${\mathfrak g}$ is a {\em ${\mathfrak
 g}$-module
 of exponential type} under the adjoint action, which means that all 
 roots of the ${\mathfrak g}$-module ${\mathfrak g}$ have the form 
 \[
 \Psi(X) = (1 + i \alpha) \lambda(X) 
 \]
 with $\lambda$ a real linear functional and $\alpha \in \RR$
 \cite[Chap. I]{BeCo}.
 It follows that the submodule ${\mathfrak n}$
 is also of exponential type.
 Passing to the dual yields that ${\mathfrak n}^*$ is a ${\mathfrak g}$-
 module of exponential type under the coadjoint action.
 Hence we obtain for the canonically induced coadjoint
 action ${\rm Ad}_G^*$ of $G$ on ${\mathfrak n}^* \cong
 {\mathfrak g}^*/{\mathfrak n}^\bot$ that all $G$-orbits in 
 ${\mathfrak n}^*$ are locally closed \cite[Chap. I, Theor\'eme 3.8]{BeCo}.
 But then the orbit space $ {\mathfrak n}^* / {\rm Ad}_G^*(G)$ is countably
 separated, by Glimm's Theorem (e.g., \cite[Chap. I, Remarque 3.9]{BeCo}).

 On the other hand, let
 \[
 \kappa: {\mathfrak n}^* / {\rm Ad}_N^*(N) \to \widehat{N}
 \]
 denote the Kirillov map, which is a homeomorphism.
 Then it is straightforward to check that $\kappa$ intertwines
 the action of ${\rm Ad}_G^*$ with the dual action, thus inducing
 a homeomorphism of orbit spaces
 \[
 {\mathfrak n}^* /{\rm Ad}_G^*(G) \to \widehat{N}/G~~.
 \]
 Hence the right-hand side is countably separated, and $N$ is regularly 
 embedded.
\end{prf}

For the formulation of the next theorem recall that the {\em nilradical}
$N$ of a solvable Lie group $G$
is defined as the maximal connected nilpotent normal subgroup of
$G$. Hence $N \lhd G$ is simply connected, with $G/N \cong
\RR^n$ \cite[Chapter III]{AuMo}.
Recall also that nilpotent, or more generally, exponential Lie groups
are of type I \cite{Ta}.
A class R solvable Lie group is defined by the requirement that for all
$x \in G$ and for all eigenvalues $\lambda$ of ${\rm Ad}(x)$, $|x|=1$
\cite{AuMo}. By contrast, exponential Lie groups are characterised by
the property that no eigenvalue of any ${\rm Ad}(x)$ is purely
imaginary \cite[Th\'eor\`eme 2.1]{BeCo}.
\begin{thm} \label{thm:wPW_solv}
 Let $G$ be a solvable Lie group,
 and let $N \lhd G$ denote the nilradical. Assume that $G$ is
 type $I$ and that $N$ is regularly embedded. Then $G$ has wPW.
 In particular, $G$ has wPW if it is exponential, or if it is of class R and
 type I.
\end{thm}
\begin{prf}
 wPW for $N$ is established by straightforward iterated application of 
 Corollary \ref{cor:abelian_ext} to a Jordan-H\"older series of $N$,
 observing that normal subgroups in nilpotent Lie groups are regularly
 embedded by Lemma
 \ref{lemma:nilp_reg_emb}.
 Moreover, $G/N \cong \RR^n$, and $N$ is type I. Hence
 Corollary \ref{cor:abelian_ext} once again applies to yield wPW for $G$.

 Now if $G$ is exponential, it is type I by 
 \cite{Ta}, and $N$ is regularly embedded by Lemma \ref{lemma:nilp_reg_emb}.
 If $G$ is of type I and class R, $N$ is regularly embedded by 
 \cite[Chapter III, Theorem 1]{AuMo}.
\end{prf}

\begin{cor}
 If $G$ is a solvable CCR Lie group, it has wPW.
\end{cor}
\begin{prf}
 CCR groups are of type I, and solvable CCR Lie groups are
 in addition of class R \cite[Chapter V, Theorem 1]{AuMo}.
 Hence the previous theorem applies.
\end{prf}

Let us next give a class of group extensions that fail to have wPW,
namely those where the normal subgroup is (nontrivial and) compact.
For this purpose, an alternative formulation of wPW, which has the
additional advantage of applying also to the non-type I setting, is
observed:
\begin{rem}
 If $G$ is type I, then the following conditions are equivalent:
\begin{enumerate}
 \item[(i)] $G$ has wPW.
 \item[(ii)] Every nonzero $f \in {\rm L}_c^\infty(G)$ is cyclic for the two-sided representation of $G$ acting on ${\rm L}^2(G)$.
 \item[(iii)] 
 For all nonzero $f \in {\rm L}_c^\infty(G)$ and nonzero every two-sided invariant operator
 $T$ on ${\rm L}^2(G)$, $Tf \not= 0$. 
\end{enumerate}
$(ii) \Leftrightarrow (iii)$ is \cite[I.I.4]{DiW}. For $(i) 
\Rightarrow (iii)$ let $f \in {\rm L}_c^\infty(G)$ and let
  $T$ denote a two-sided invariant operator. Under the Plancherel transform,
 \[
 T \simeq \int_{\widehat{G}}^\oplus m(\sigma)  \cdot {\rm Id}_{
 {\mathcal H}_{\sigma} \otimes \overline{\mathcal H}_\sigma} ~d\nu_G(\sigma)
 \]
 for a certain Borel mapping $m \in {\rm L}^\infty(\widehat{G})$.
 If $T \not= 0$, $m$ does not vanish identically.
 But then $(i)$ implies that
 \[ (Tf)^\wedge (\sigma) = m(\sigma) \widehat{f}(\sigma) \]
 does not vanish identically, thus $Tf \not=0$.
 For $(iii) \Rightarrow (i)$ assume that $\widehat{f}$
 vanishes on a set $\Sigma \subset \widehat{N}$ of 
 positive Plancherel measure. Let $P$ denote the projection 
 defined by 
 \[ 
 P \simeq \int_{\widehat{G}} \chi_{\Sigma}(\sigma) \cdot
 {\rm Id}_{{\mathcal H}_{\sigma} \otimes \overline{\mathcal H}_\sigma}~~.
 \]  
 Then $P$ is two-sided invariant, nontrivial, but $Pf=0$.

 Similar arguments apply to show that the following conditions 
 are equivalent, for regularly embedded $N \lhd G$:
\begin{enumerate}
 \item[(i)] Condition (\ref{eqn:inv_wPW}) holds.
 \item[(ii)] Every nonzero $f \in {\rm L}_c^\infty(N)$ is cyclic for the von-Neumann 
 algebra generated by the two-sided representation of $N$ and the 
 representation of $G$ acting on ${\rm L}^2(N)$ by conjugation.
 \item[(iii)] 
 For all nonzero $f \in {\rm L}_c^\infty(N)$ and nonzero two-sided
 invariant operator $T$ on ${\rm L}^2(N)$ commuting with the conjugation
 action of $G$, $Tf \not= 0$. 
\end{enumerate}
\end{rem}

\begin{prop}
\label{prop:cp_norm}
 If $G$ has wPW, it has no nontrivial compact normal subgroups.
\end{prop}
\begin{prf}
 Assume that $K \lhd G$ is compact, and
 consider the subspace 
 \begin{eqnarray*}
 {\rm L}^2_K(G) & = &  \{ f \in {\rm L}^2(G) : ~\forall k \in K:  f(xk) = f(x)
 \} \\
   & = & \{ f \in {\rm L}^2(G) : ~\forall k \in K:  f(kx) = f(x)
 \} ~~. 
\end{eqnarray*}
 It is easy to see that ${\rm L}^2_K(G)$ is closed and two-sided invariant.
 Moreover clearly ${\rm L}^2_K(G) \not= {\rm L}^2(G)$ and 
 ${\rm L}^2_K(G) \cap {\rm L}_c^\infty \not= \{ 0 \}$. Hence
 $G$ does not have wPW, by the previous remark. 
\end{prf}

Proposition \ref{prop:cp_norm} in fact follows from \cite[Lemma 1.1]{KaLaSc},
where it is stated for the topological Paley-Wiener condition. Note that
the topological Paley-Wiener condition is weaker than wPW.
\cite[Theorem 3.2]{KaLaSc} shows that for SIN groups, the two conditions 
coincide with the necessary condition derived in the
previous proposition.


\begin{ex}
 An example where condition (\ref{eqn:inv_wPW}) holds, but $N$ does not have wPW is constructed
as follows: Consider $G = \QQ_p \semdir \QQ_p^\times$, where $\QQ_p$ denotes the
field of $p$-adic numbers, and its unit group $\QQ_p^\times$
acts by multiplication. $\QQ_p$ is self-dual, and the dual action of $\QQ_p^\times$ 
is again by multiplication. In particular, $\widehat{\QQ_p}$ consists of the two
dual orbits $\{ 0 \}$ (which has measure zero) and $\QQ_p^\times$. Moreover,
the action of $\QQ_p^\times$ on the large orbit is free. Hence
Kleppner and Lipsman's theorem yields that the Plancherel measure is supported on
a single point, and the wPW property is an immediate consequence of the Plancherel
theorem. (Of course, Theorem \ref{thm:main} also applies, with both conditions
trivially fulfilled.)

On the other hand, $\QQ_p$ has the nontrivial compact subgroup $\ZZ_p$,
hence $\QQ_p$ does not have wPW.
%
\end{ex}

\section{Acknowledgements} I thank G\"unter Schlichting for interesting
discussions and for a copy of \cite{KaLaSc}.

\bibliography{wpw.bib}

\begin{thebibliography}{10}

\bibitem{ArnLud}
D.~Arnal and J.~Ludwig.
\newblock {Q}.{U}.{P}. and {P}aley-{W}iener properties of unimodular,
  especially nilpotent, {L}ie groups.
\newblock {\em Proc. Am. Math. Soc.}, 125:1071--1080, 1997.

\bibitem{AuMo}
L.~Auslander and C.C. Moore.
\newblock {\em Unitary representations of solvable {L}ie groups}, volume~62 of
  {\em Mem. Am. Math. Soc.}
\newblock 1966.

\bibitem{BaKl}
L.~Baggett and A.~Kleppner.
\newblock Multiplier representations of abelian groups.
\newblock {\em J. Funct. Anal.}, 14:299--324, 1973.

\bibitem{BeCo}
P.~Bernat, N.~Conze, M.~Duflo, M.~L\'evy-Nahas, M.~M.~Ra\"{\i}s, P.~Renouard,
  and M.~Vergne.
\newblock {\em Repr\'esentations des Groupes de {L}ie R\'esolubles}.
\newblock Dunod, Paris, 1972.

\bibitem{DiW}
J.~Dixmier.
\newblock {\em Von {N}eumann-Algebras}.
\newblock North-Holland, 1981.

\bibitem{DuMo}
M.~Duflo and C.C. Moore.
\newblock On the regular representation of a nonunimodular locally compact
  group.
\newblock {\em J. Functional Analysis}, 21:209--243, 1976.

\bibitem{EchKanKum}
S.~Echterhoff, E.~Kaniuth, and A.~Kumar.
\newblock A qualitative uncertainty principle for locally compact groups.
\newblock {\em Forum Math.}, 3:355--369, 1991.

\bibitem{Fo}
G.B. Folland.
\newblock {\em A course in abstract harmonic analysis}.
\newblock Studies in advanced mathematics. CRC Press, 1995.

\bibitem{Gar1}
G.~Garimella.
\newblock Un th\'eor\`eme de {P}aley-{W}iener pour les groupes de {L}ie
  nilpotents.
\newblock {\em J. {L}ie {T}heory}, 5:165--172, 1995.

\bibitem{Gar2}
G.~Garimella.
\newblock Weak {P}aley-{W}iener property for completely solvable {L}ie groups.
\newblock {\em Pac. J. Math.}, 187:51--63, 1999.

\bibitem{Hog}
J.A. Hogan.
\newblock A qualitative uncertainty principle for locally compact abelian
  groups.
\newblock {\em Proc. Centre for Math. Analysis}, 16:133--142, 1988.

\bibitem{KaLaSc}
E.~Kaniuth, A.T. Lau, and G~Schlichting.
\newblock A topological {P}aley-{W}iener property for locally compact groups.
\newblock Preprint.

\bibitem{KlLi}
A.~Kleppner and R.L. Lipsman.
\newblock The {P}lancherel formula for group extensions.
\newblock {\em Ann. Sci. Ec. Norm. Sup.}, 4:459--516, 1972.

\bibitem{Li}
R.L. Lipsman.
\newblock {\em Group representations}, volume 388 of {\em Lecture Notes in
  Mathematics}.
\newblock Springer, 1974.

\bibitem{LipRos}
R.L. Lipsman and J.~Rosenberg.
\newblock The behavior of {F}ourier transforms for nilpotent {L}ie groups.
\newblock {\em Trans. Amer. Math. Soc.}, 348:1031--1050, 1996.

\bibitem{Ma}
G.W. Mackey.
\newblock Induced representations of locally compact groups {I}.
\newblock {\em Ann. of Math.}, 55:101--139, 1952.

\bibitem{Mo}
J.D. Moss.
\newblock A {P}aley-{W}iener theorem for selected nilpotent {L}ie groups.
\newblock {\em J. Funct. Anal.}, 114:395--411, 1993.

\bibitem{Pa}
R.~Park.
\newblock A {P}aley-{W}iener theorem for all two- and three-step nilpotent
  {L}ie groups.
\newblock {\em J. Funct. Anal.}, 133:277--300, 1995.

\bibitem{Ta}
O.~Takenouchi.
\newblock Sur la facteur-repr\'esentation d'un groupe de {L}ie r\'esoluble de
  type ({E}).
\newblock {\em Math. J. Okayama Univ.}, 7:151--161, 1957.

\end{thebibliography}
\bibliographystyle{plain}

\end{document}